\documentclass[12pt]{article}
\usepackage{amsfonts}
\newtheorem{theorem}{Theorem}[section]
\newtheorem{lemma}[theorem]{Lemma}
\newtheorem{prop}[theorem]{Proposition}

\newtheorem{unnumber}{}

\newtheorem{prepf}[unnumber]{Proof}
\newenvironment{pf}{\prepf\rm}{\endprepf}
\newtheorem{preprob}[unnumber]{Problem}
\newenvironment{prob}{\preprob\rm}{\endpreprob}

\newcommand{\qed}{\hfill$\Box$}

\begin{document}

\title{Groups with right-invariant multiorders}
\author{Peter J. Cameron\\School of Mathematical Sciences\\
Queen Mary, University of London\\Mile End Road\\London E1 4NS, UK}
\date{}
\maketitle

\begin{abstract}
A \emph{Cayley object} for a group $G$ is a structure on which $G$ acts
regularly as a group of automorphisms. The main theorem asserts that a
necessary and sufficient condition for the free abelian group $G$ of rank $m$
to have the generic $n$-tuple of linear orders as a Cayley object is that
$m>n$. The background to this theorem is discussed. The proof uses Kronecker's
Theorem on diophantine approximation.
\end{abstract}

\section{Cayley objects and homogeneous structures}

The \emph{regular representation} of a group $G$ is the representation of
the group acting on itself by right multiplication. A \emph{Cayley object}
for $G$ is a structure on $G$ admitting the regular representation as a
group of automorphisms. The name comes from the fact that a Cayley graph
for $G$ is precisely a Cayley object which happens to be a graph.

A Cayley object must admit a transitive automorphism group. There is some
interest in investigating objects with a high degree of symmetry which are
Cayley objects for a group, or (in the other direction) groups which have
a given highly symmetric object as a Cayley object. This is the topic
of~\cite{hco}; I refer to that paper for further motivation.

All objects here will be relational structures, consisting of a set carrying
a collection of relations of various arities. A substructure of a relational
structure will always be the induced substructure on a subset, consisting of
all instances of each relation such that all arguments lie within the subset.

A relational structure $M$ is said to be \emph{homogeneous} if any isomorphism
between finite substructures can be extended to an automorphism of $M$. This
will be our ``strong symmetry condition''.

The \emph{age} of a relational structure $M$ is the class of all finite
relational structures of the same type which can be embedded into $M$.

Fra\"{\i}ss\'e~\cite{fr} gave a necessary and sufficient condition for a 
class $\mathcal{C}$ of finite structures to be the age of a countable
homogeneous structure:
\begin{enumerate}
\item $\mathcal{C}$ is closed under isomorphism;
\item $\mathcal{C}$ is closed under taking substructures;
\item $\mathcal{C}$ contains only finitely many members up to isomorphism;
\item $\mathcal{C}$ has the \emph{amalgamation property}, that is, given
$A,B_1,B_2\in\mathcal{C}$ with embeddings $f_i:A\to B_i$ for $i=1,2$, there
exists $C\in\mathcal{C}$ and embeddings $g_i:B_i\to C$ for $i=1,2$ such
that the composite embeddings $g_1f_1$ and $g_2f_2$ agree.
\end{enumerate}
Moreover, if these conditions hold, there is a unique countable homogeneous
structure $M$ with age $\mathcal{C}$ (up to isomorphism). Such a class
$\mathcal{C}$ is called a \emph{Fra\"{\i}ss\'e class}, and $M$ is its
\emph{Fra\"{\i}ss\'e limit}.

We say that $\mathcal{C}$ has the \emph{strong amalgamation property} if the
amalgamation can be done without identifying points outside $A$: that is, if
$b_1\in B_2$ and $b_2\in B_2$ satisfy $g_1(b_1)=g_2(b_2)$, then there
exists $a\in A$ such that $b_i=f_i(a)$ for $i=1,2$.

For example, the class of all finite totally ordered sets is a Fra\"{\i}ss\'e
class; its Fra\"{\i}ss\'e limit is the ordered set $\mathbb{Q}$, the unique
countable dense ordered set without endpoints. I generalise this example in
the next section.

The homogeneous structure $M$ with age $\mathcal{C}$ is characterised by the
following \emph{extension property}:
\begin{quote}
If $A,B\in\mathcal{C}$ with $A\subseteq B$ and $|B|=|A|+1$, then every
embedding of $A$ into $M$ can be extended to an embedding of $B$ into $M$.
\end{quote}

\section{Multiorders}

An \emph{$n$-order} is a set with $n$ linear orders. If we
do not need to specify the number of orders, we refer to a \emph{multiorder}.

The class of finite $n$-orders is a Fra\"{\i}ss\'e class. (More generally,
if we take any finite number of Fra\"ss\'e classes, each of which has strong
amalgamation, and consider the finite sets carrying a structure from each
class, with no relationship between the different structures, we obtain a
Fra\"{\i}ss\'e class.)

The Fra\"{\i}ss\'e limit of the class of finite $n$-orders will be called
the \emph{generic} (countable) $n$-order.

The case $n=2$ arises in connection with the thriving field of
\textit{permutation patterns}. If a finite
set $X$ carries a $2$-order, we can use the first order to enumerate $X$
as $(x_1<_1x_2<_1\cdots<_1x_k)$, and then the second order defines a
permutation of the labels $\{1,2,\ldots,n\}$. The notion of induced
substructure coincides exactly with that used in the theory of permutation
patterns. So, in a sense, the theory of permutation patterns is the theory of
the age of the generic countable $2$-order. Is there a similar theory for
the generic $n$-orders with $n>2$?

In this context, the countable homogeneous $2$-orders were determined
in~\cite{hp}.

\begin{prob}
Determine all countable homogeneous $n$-orders, for $n>2$.
\end{prob}

The extension property characterising the generic $n$-order on a countable
set $X$ is the following:
\begin{quote}
Given any $k$ points $x_1,\ldots,x_k$ of $X$, for each $i\in\{1,\ldots,n\}$,
let $I_i$ be one of the $k+1$ intervals (including semi-infinite intervals)
into which $X$ is divided by $x_1,\ldots,x_k$ in the order $<_i$. Then
$I_1\cap\cdots\cap I_n\ne\emptyset$.
\end{quote}
This is because adding a point to a finite totally ordered set involves
putting it into one of the intervals defined by the set: i.e. before the
first element, or between the $i$th and $(i+1)$st for $i=1,\ldots,k-1$, or
after the last element.

For ease of use, we give a simpler but equivalent condition.

\begin{prop}
An $n$-order $(<_1,\ldots,<_n)$ on a countable set $X$ is generic if and only
if, for any choice of $x_i$ and $y_i$ (for $i=1,\ldots,n$) with
$x_i<_iy_i$ (possibly $x_i=-\infty$ or $y_i=\infty$), there is a point
$z\in X$ satisfying $x_i<_iz<_iy_i$ for $i=1,\ldots,n$.
\end{prop}

\begin{pf}
It is clear that the condition in the Proposition implies that in the 
extension property. Suppose that the condition in the extension property
is true, and assume the hypotheses of the Proposition. Take the finite
set $\{x_1,\ldots,x_n,y_1,\ldots,y_n\}$. For each $i$, the set of points
$z$ satisfying $x_i<z<y_i$ is a union of intervals defined by this finite
set; pick one of them. By the extension property, the intersection of the
chosen intervals is non-empty.\qed
\end{pf}

We will be interested in the case where all the orders are right-invariant
for a countable group. This means that none of them have end-points, and
we don't need to worry about the ``semi-infinite'' intervals.

\section{Dense right orders on groups}

To say that a group $G$ has a totally ordered set which is a Cayley object
means that there is a total order $<$ on $G$ which is \emph{right-invariant},
that is, if $x<y$, then $xg<yg$ for any $g\in G$. If we have such an order,
let $P=\{g\in G:1<g\}$; then
\begin{itemize}
\item[(a)] $G$ is the disjoint union of $\{1\}$, $P$ and $P^{-1}$;
\item[(b)] $P^2\subseteq P$.
\end{itemize}
Conversely, if we have a set $P$ satisfying these two conditions then,
setting $x<y$ if $y=px$ for some $p\in P$ defines a right-invariant order
on $G$. Moreover, the order is dense if and only if (b) is replaced by the
stronger condition
\begin{itemize}
\item[(bb)] $P^2=P$.
\end{itemize}
For if $x<y$, then $x=py$ for some $p\in P$. If $P=P^2$, then write $p=qr$
for some $q,r\in P$; then $x<rx<qrx=y$.

A group is said to be \emph{right-orderable} if it has a right-invariant order
(sometimes called a \emph{right order} for short). A great deal is known about
right-orderable groups (see Chapter VII of~\cite{mr} for a survey, and note
that since a right order of an abelian group is also a left order, the results
of the whole book apply in the case of abelian groups). Less attention has
been paid to groups with a dense right order. Here is the example which will
be important to us.

\begin{theorem}
Let $\mathbb{Z}^m$ denote the free abelian group of rank $m>1$. Suppose that
$<$ is a right order on $G$. Then there is a non-zero vector $c\in\mathbb{R}^m$
such that $x<y$ if $c.x<c.y$, where the dot denotes the usual
inner product. Moreover, if the components of $c$ are linearly independent
over $\mathbb{Q}$, then the order is dense, and $x<y$ if and only if
$c.x<c.y$.
\end{theorem}

Note that, if the components of $c$ are not linearly independent over 
$\mathbb{Q}$, then there are non-zero elements $z$ of $\mathbb{Z}^m$ which
satisfy $c.z=0$, forming a subgroup $A$ which is free abelian of smaller rank;
to complete the specification of the order, we have to choose a right order
of $A$. Note that the order is \emph{non-archimedean} in this case; if $a,b$
are positive elements with $a\in A$ and $b\notin A$, then $a^n<b$ for all
positive $n$.

For example, $\mathbb{Z}$ has just two right orders (the usual order and its
reverse), neither of which is dense. For $\mathbb{Z}^2$, using a vector
of the form $c=(1,\alpha)$ gives a dense order if $\alpha$ is irrational.
However, if $\alpha$ is rational, or if $c=(0,1)$, then we do not yet have
enough information to define the order, since the points $z\in\mathbb{Z}^2$
which satisfy $c.z=0$ will form a subgroup whose order is not yet specified.
This subgroup has rank $1$, and so (as before) has just two orders.

I have not found a convenient exposition of the proof of this theorem, so here
is a sketch. By factoring out the subgroup of ``small'' elements, we may
assume that the ordering is archimedean. Then a theorem of H\"older~\cite{ho}
shows that there is an isomorphism to an additive subgroup of $\mathbb{R}$,
which clearly has the form given. See also~\cite[Theorem 1.3.4]{mr} or
\cite[p. 62]{fs}.

\section{The main theorem}

The main result of this paper is the first known class of groups admitting
homogeneous right multiorders. This result was conjectured in \cite{hco}.

\begin{theorem}
Let $m$ and $n$ be positive integers. The free abelian group $\mathbb{Z}^m$ of
rank $m$ has a right-invariant generic $n$-tuple of orders if and only if
$m>n$.
\end{theorem}

The proof of the theorem requires a number of lemmas. First we show that,
if $m>n$, then there
is a $\mathbb{Z}^m$-invariant generic $n$-tuple of orders. We note first that
it suffices to show the result when $m=n+1$, since dropping some orders from
a $G$-invariant generic multiorder yields a $G$-invariant generic multiorder.

The proof uses an important result of Kronecker~\cite{kro} on diophantine
approximation, for which several proofs are given in Chapter XXIII of Hardy
and Wright~\cite{hw}.

\begin{theorem}
Let $m$ be a positive integer, and let $c\in\mathbb{R}^m$ be a vector whose
components are linearly independent on $\mathbb{Q}$. Then, given any
$\epsilon>0$, any line in $\mathbb{R}^m$ with direction vector $c$ passes
within distance $\epsilon$ of some lattice point in $\mathbb{Z}^m$.
\end{theorem}

We also need an existence result for a certain kind of matrix.

\begin{lemma}
Let $m$ be a positive integer. Then there exists a $m\times m$ real matrix
$A$ having the properties
\begin{enumerate}
\item $A$ is invertible;
\item each row of $A$ has components which are linearly independent over
$\mathbb{Q}$;
\item the last row of $A$ is orthogonal to all the others.
\end{enumerate}
\label{mx}
\end{lemma}

\begin{pf}
The set of $n\times n$ invertible matrices is a complete metric space. (A
little care is required; we take $n^2+1$ coordinates, the matrix entries and
the inverse of the determinant, the latter required to ensure that a Cauchy
sequence of invertible matrices cannot converge to a singular matrix.) Then
the condition (c) defines a closed subspace, which is therefore also complete.

Now condition (b) restricts us to a countable intersection of subsets (one
for each choice of row of the matrix and coefficients in a rational linear
combination of the entries in that row). Each such set is obviously open
(since we require that the linear combination is not zero). We claim that
each such set is dense. For take a matrix $A$ satisfying conditions (a) and
(c), such that $\sum q_ja_{ij}=0$. Suppose, without loss, that $a_{ij}\ne0$.
Take $t$ to be a real number which is transcendental over the field generated
by all the matrix entries and is arbitrarily close to $1$. Replace $a_{1k}$ 
by $ta_{1k}$ for $k=1,\ldots,n-1$, and $a_{1n}$ by $t^{-1}a_{1n}$. Then
condition (c) is preserved; (a) is preserved if $t$ is sufficiently close
to $1$; and the linear combination is no longer zero.

The set of matrices satisfying all three conditions is a residual set (a
countable intersection of open dense sets) in a complete metric space, and
hence is non-empty, by the Baire Category Theorem (\cite{ba}; see~\cite{ox}
for discussion).\qed
\end{pf}

\paragraph{Remark} It may not be too difficult to write down explicit examples
of matrices with these properties. For example, when $n=2$, we can take
\[A=\pmatrix{1&\alpha\cr-\alpha&1\cr}\]
for any irrational number $\alpha$.

\medskip

Now we give the construction. Let $A$ be a matrix having the properties of
Lemma~\ref{mx}. Use the first $m-1$ rows to define $m-1$ right-invariant
total orders $<_1,\ldots,<_{m-1}$ on $\mathbb{Z}^m$. We claim that this
$(m-1)$-tuple is generic.

An interval in the $i$th order consists of the vectors lying between two
parallel hyperplanes perpendicular to the $i$th row of the matrix. Since the
matrix is invertible, the intersection of $m-1$ intervals (one for each
order is a cylinder with parallelepiped cross-section in a direction 
orthogonal to the first $m-1$ rows of the matrix, hence (by condition (c))
parallel to the $m$th row. By Kronecker's Theorem, there is a lattice point
arbitrarily close to this line, and in particular close enough that it lies
in the cylinder defined by the intervals. So this intersection is non-empty
in the lattice $\mathbb{Z}^m$, and we are done.

\medskip

Now we turn to the non-existence proofs.

Since we may omit some orders from a generic multiorder and it remains
generic, we may assume that $m=n$. Our proof is by induction on $n$; it
is split into two cases, of which only the second case requires the
induction hypothesis.

Let $(<_1,\ldots,<_n)$ be an $n$-tuple of right-invariant linear orders 
on $\mathbb{Z}^n$. We have to prove that this tuple is not generic. Let
$c_1,c_2,\ldots,c_n$ be vectors defining the top section of the ordering.
We use this notation for the remainder of the proof.

\begin{lemma}
If $c_1,\ldots,c_n$ are linearly dependent, then the $n$-tuple of orders is
not generic.
\end{lemma}

\begin{pf} Suppose that $c_k$ is a linear combination of $c_1,\ldots,c_{k-1}$,
say $c_k=a_1c_1+\cdots+a_{k-1}c_{k-1}$. By reversing some of the orders if
necessary, we may assume that all the coefficients are non-negative. Now
choose intervals $x_i\le c_i.z\le y_i$ in the group. Any point $z$ lying in
all these intervals must also lie in the interval
\[\sum_{j=1}^{k-1}a_jx_j\le c_k.z\le\sum_{j=1}^{k-1}a_jy_j.\]
So the interval $x_k\le c_k.z\le y_k$ does not meet the intersection of these
$k-1$ intervals if we choose, say, $y_k<\sum_{j=1}^{k-1}a_jx_j$.\qed
\end{pf}

\begin{lemma}
If $c_1,\ldots,c_k$ are linearly independent and at least one of them has
linearly dependent components over $\mathbb{Q}$, then the $n$-tuple of
orders is not generic.
\end{lemma}

\begin{pf}
Without loss of generality, we may assume that $c_1$ has linearly dependent
components over $\mathbb{Q}$. Then $A=\{z\in\mathbb{Z}^n:c_1.z=0\}$ is a
non-zero subgroup of $\mathbb{Z}^m$, and contains an interval $I_1$ in the
order $<_1$. So the restrictions of the other orders to $A$ form an
$(n-1)$-tuple of orders on an abelian group of rank at most $n-1$. By the
inductive hypothesis, they cannot be generic, so some intersection of
intervals in these orders is disjoint from $A$, and hence from $I_1$. So
the original order is not generic.\qed
\end{pf}

\begin{lemma}
If $c_1,\ldots,c_k$ are linearly independent and all of them have components
which are linearly dependent over $\mathbb{Q}$, then the $n$-tuple of orders
is not generic.
\end{lemma}

\begin{pf}
Each of the orders $<_1,\ldots,<_n$ is dense; an interval in $<_i$ consists
of the lattice points lying between two parallel hyperplanes perpendicular
to $c_i$, and these hyperplanes may be arbitrarily close together. So the
intersection of the $n$ intervals is a parallelepiped whose volume can be
made arbitrarily small. This parallelepiped tiles the Euclidean space, so
if we make its volume less than $1$ we can find a translate containing no
lattice point.\qed
\end{pf}

These three lemmas complete the proof of the theorem.\qed

\begin{prob}
Find further examples of groups with generic right multiorders.
\end{prob}


\begin{thebibliography}{99}

\bibitem{ba}
R. Baire, Sur les fonctions de variables r\'eeles,
\textit{Ann. di Mat.} \textbf{3} (1899), 1--123.

\bibitem{mr}
R. T. Botto Mura and A. H. Rhemtulla, \textit{Orderable Groups},
Marcel Dekker, New York, 1977.

\bibitem{hco}
P. J. Cameron, Homogeneous Cayley objects, 
\textit{European J. Combinatorics} \textbf{21} (2000), 745--760. 

\bibitem{hp}
P. J. Cameron, Homogeneous permutations,
\textit{Electronic J. Combinatorics} \textbf{9(2)} (2002), \#R2 (9pp).

\bibitem{fr}
R. Fra\"{\i}ss\'e,
Sur certains relations qui g\'en\'eralisent l'ordre des nombres rationnels,
\textit{C. R. Acad. Sci. Paris} \textbf{237} (1953), 540--542.

\bibitem{fs}
L. Fuchs and L. Salce, \textit{Modules over non-Noetherian domains},
American Math. Soc., Providence, R.I., 2001.

\bibitem{hw}
G. H. Hardy and E. M. Wright,
\textit{An Introduction to the Theory of Numbers}, fifth edition, Oxford
University Press, Oxford, 1979.

\bibitem{ho}
O. H\"older, 
Die Axiome der Quantit\"at und die Lehre vom Mass,
\textit{Berichte Verh. s\"achs. Gesell. Wiss.} (Leipzig), Math. Phys. Cl.
\textbf{53} (1901), 1--64.

\bibitem{kro}
L. Kronecker, Die Periodensysteme von Funktionen reeller Variablen,
\textit{Berliner Sitzungsberichte} (1884), 1071--1080.

\bibitem{ox}
J. C. Oxtoby, \textit{Measure and Category},
Springer, Berlin, 1980.

\end{thebibliography}
\end{document}